\documentclass[12pt]{article}
\usepackage{amsmath}
\usepackage{amsfonts}
\usepackage{amssymb}
\usepackage{latexsym}

\baselineskip
\newtheorem{thm}{Theorem}
\newtheorem{cor}{Corollary}
\newtheorem{defn}{Definition}
\newtheorem{lem}{Lemma}

\input{tcilatex}

\begin{document}

\title{(U,V)-ordering and a duality theorem for risk aversion and
Lorenz-type orderings}
\author{A. Giovagnoli \\
University of Bologna, Italy \and H.P. Wynn \\
London School of Economics, UK }
\maketitle

\begin{abstract}
There is a duality theory connecting certain stochastic orderings between
cumulative distribution functions $F_{1},F_{2}$ and stochastic orderings
between their inverses $F_{1}^{-1},F_{2}^{-1}$. This underlies some theories
of utility in the case of the cdf and deprivation indices in the case of the
inverse. Under certain conditions there is an equivalence between the two
theories. An example is the equivalence between second order stochastic
dominance and the Lorenz ordering. This duality is generalised to include
the case where there is \textquotedblleft distortion\textquotedblright\ of
the cdf of the form $v(F)$ and also of the inverse. A comprehensive duality
theorem is presented in a form which includes the distortions and links the
duality to the parallel theories of risk and deprivation indices. It is
shown that some well-known examples are special cases of the results,
including some from the Yaari social welfare theory and the theory of
majorization.
\end{abstract}

\noindent \textbf{Key Words\qquad }Income inequality, Prospect Theory,
Stochastic orderings, Utility theory, Yaari's Functionals

\noindent \textbf{Subject Codes \ }C020, D690; D390\bigskip

\section{Introduction}

Many results in mathematical utility theory and in the parallel theories of
poverty can be cast in terms of stochastic orderings; a standard reference
is Shaked and Shanthikumar (2007). Moreover, it is becoming clear from work
in these areas that there is a duality between orderings based on income
distributions expressed via the cumulative income distribution (cdf) $F$ and
certain orderings on the quantile function $F^{-1}$; see, in particular,
Yaari (1987). Such matters are also of growing interest in financial and
insurance risk areas. A useful text is M\"{u}ller and Stoyan (2002). An
example of the duality is that between second order stochastic dominance
(SSD) and so-called weak Lorenz ordering: 
\begin{equation}
\int_{-\infty }^{x}F_{1}dt\geq \int_{-\infty }^{x}F_{2}dt\;\;\mbox{for all}%
\;x\in \mathbf{R}\Leftrightarrow \int_{0}^{\alpha }F_{1}^{-1}dt\leq
\int_{0}^{\alpha }F_{2}^{-1}dt\;\;\mbox{for all}\;0\leq \alpha \leq 1.
\label{1}
\end{equation}%
This result has been studied in important contributions by Muliere and
Scarsini (1999) and Ogryczak and Ruszczynski (2002), who show the link to
Fenchel duality. The book by Pallaschke and Rolewicz (1997) covers relevant
generalisations of Fenchel duality. Sordo and Ramos (2007) cover similar
ground with a discussion of the literature, going back to Lorenz (1905). A
generalisation of (\ref{1}) is at the heart of this paper.

There is also an equivalence between orderings and their characterizations
in terms of a set of order preserving functions, expected utilities in the
language of economics. Thus, first order stochastic dominance defined as $%
F_{1}(x)\geq F_{2}(x)$ for all $x$, holds if and only if 
\begin{equation*}
E_{X}u(X)\leq E_{Y}u(Y)
\end{equation*}%
where $X\sim F_{1}$ and $Y\sim F_{2},$ for all non-decreasing $u(\cdot )$.
For second order stochastic dominance (SSD) the equivalence is for all
non-decreasing concave $u(\cdot )$ for which expectations exist. There is a
utility theory associated with $F^{-1}$ which is often referred to as the
dual theory of \textit{rank-dependent} utilities, and studies concepts such
as \textquotedblleft deprivation indexes\textquotedblright : see\ Atkinson
(1970), Sen (1973), Cowell (1977). For the moment we simply note that these
dual utilities measure the rank position of an individual (in terms of
income) in the population, as opposed to the actual level of income. There
are a variety of duality \emph{theorems} relating the utility theories based
on $F$ and $F^{-1}$. Recent examples are Maccheroni, Muliere and Zoli (2005)
and Chateauneuf and Moyes (2005).

A rich and related area is the study of \textquotedblleft
distortions\textquotedblright\ on a cumulative distribution function $F$. It
is fundamental to our approach that there are two types of distortion and
that they can be applied to $F$ or to $F^{-1}$. One type takes the form of a
direct transformation of the cdf: $v(F)$, for some function $v(\cdot )$. The
other type is applied to the base measure so that, for an integrable
function $g(\cdot )$, $\int g(t)dt$ becomes $\int g(t)du(t)$. For $F^{-1}$
we have analogous distortions. Distortions have been studied in Prospect
Theory using terminology such as \textquotedblleft probability
weighting\textquotedblright , see Kahneman and Tversky (1979), (1992).

It is a main aim of the paper to present a general type of stochastic
ordering which helps to unify the above theories. This is done via our main
duality theorem, Theorem 1 in Section 3. The proof, in Appendix 1, uses the
promised generalisation of (\ref{1}), Lemma 1, which incorporates
distortions of both types just mentioned, applied both to $F$ and, swapped
over, to $F^{-1}$. Furthermore, each side of the duality has an equivalent
representation in terms of utility functions which, perhaps surprisingly,
uses the same class of functions as are used for the distortion. This means
that Theorem 1 has four equivalent parts. Our stochastic ordering is defined
given two base (distortion) functions $u_{0}$ and $v_{0}$. The duality
involving $F$ and $F^{-1}$ and the utility versions are then fixed. The
proofs are somewhat technical because we assume general cdf's. Appropriate
forms of some standard results are needed, such as integration by parts, and
these are put into Appendix 2.

The test of a theory may be the range of its special cases. In Section 6 we
cover such cases and also issues such as dominated risk aversion and
inequality aversion.

\section{Orderings and duality}

We start with a general definition which involves a simultaneous distortion
of the cdf's and the measure. Let $U,V$ be classes of functions $u,v$, where 
$u:\boldsymbol{R}\rightarrow \boldsymbol{R}$ and $v:[0,1]\rightarrow \lbrack
0,1]$. All such functions will be of bounded variation, so that we can take
the associated measures. The notation $\int f(x)du(x)$ and $\int g(\alpha
)dv(\alpha )$ will be used for the Lebesgue-Stieltjes integrals, but when
the integrating variable is clear, it will be shortened to $\int f(x)du$ and 
$\int g(\alpha )dv.$ Throughout the paper it will be understood that when $u$
or $v$ is right continuous the associated Stieltjes integral will be
extended to half-open intervals inclusive of the upper end-point, whereas if
the function is left continuous the Stieltjes integral will be extended to
intervals half-open to the right.

\begin{defn}
Let $U,V$ be classes of functions $u,v$ and let $F_{1},F_{2}$ be two cdf's.
We say that $F_{1}$\textbf{\ is less than }$F_{2}$\textbf{\ in the }$(U,V)$%
\textbf{-ordering} if 
\begin{equation}
\int_{-\infty }^{\infty }v(F_{1}(x))du(x)\geq \int_{-\infty }^{\infty
}v(F_{2}(x))du(x)  \label{2}
\end{equation}%
for all $u\in U,$ $v\in V$ for which the integrals exist.\bigskip 
\end{defn}

This paper is also concerned with dual orderings which take the form given
by the following definition. For a cdf $F$ we define $F^{-1}$ in the usual
way: 
\begin{equation*}
F^{-1}(\alpha )=\inf \{x:F(x)\geq \alpha \},\text{ }\alpha \in (0,1)
\end{equation*}%
and, following the standard, take $F$ to be right continuous so that $F^{-1}$
is left continuous.

\begin{defn}
Let $\tilde{U},\tilde{V}$ be classes of functions $\tilde{u},\tilde{v}$ and
let $F_{1},F_{2}$ be two cdf's. We say that $F_{1}$\textbf{\ is less than }$%
F_{2}$\textbf{\ in the dual }$(\tilde{U},\tilde{V})$\textbf{-ordering }if 
\begin{equation}
\int_{0}^{1}\tilde{u}(F_{1}^{-1}(\alpha ))d\tilde{v}(\alpha )\leq
\int_{0}^{1}\tilde{u}(F_{2}^{-1}(\alpha ))d\tilde{v}(\alpha )  \label{3}
\end{equation}%
for all $\tilde{u}$ in $\tilde{U}$ and $\tilde{v}$ in $\tilde{V}$.\bigskip 
\end{defn}

A duality theorem in the context of Definitions 1 and 2 is a collection $%
\{U,V,\tilde{U},\tilde{V}\}$ such that the $(U,V)$-ordering and the dual $(%
\tilde{U},\tilde{V})$-ordering are mathematically equivalent. A well-known
example is for the ordinary (first order) stochastic dominance\smallskip 
\begin{equation*}
F_{1}(x)\geq F_{2}(x)\text{ for all }x\in R\iff F_{1}^{-1}(\alpha )\leq
F_{2}^{-1}(\alpha )\text{ for all }\alpha \in (0,1)
\end{equation*}%
\smallskip where $U=\{$all indicator functions ${\mathbb{I}}_{[c,\infty
)}(x)\},$ $V$ = \{identity on $[0,1]\},$ $\tilde{U}=\{$identity on $%
\boldsymbol{R}\}$, $\tilde{V}=\{$all indicator functions ${\mathbb{I}}%
_{[p,1)}(\alpha )\}.$

Importantly, using integration by parts under suitable conditions, each of
the inequalities in Definitions 1 and 2 may have an equivalent version in
terms of expected utility. In that case the $(U,V)$-ordering is equivalent
to the statement 
\begin{equation}
\int_{-\infty }^{\infty }u(x)dv(F_{1}(x))\leq \int_{-\infty }^{\infty
}u(x)dv(F_{2}(x))  \label{4}
\end{equation}%
for all $u$ in $U$ and $v$ in $V$. The dual $(\tilde{U},\tilde{V})$-ordering
is equivalent to%
\begin{equation}
\int_{0}^{1}\tilde{v}(\alpha )d\tilde{u}(F_{1}^{-1}(\alpha ))\geq
\int_{0}^{1}\tilde{v}(\alpha )d\tilde{u}(F_{2}^{-1}(\alpha ))  \label{5)}
\end{equation}%
for all $u$ in $\tilde{U}$ and $v$ in $\tilde{V}$. It would not be too
presumptuous to say that the majority of stochastic orderings defined in the
literature (see the list at the end of Shaked and Shanthikumar, 2007) are of
the form (\ref{2}), (\ref{3}), (\ref{4}) or (\ref{5)}). An equivalence
theorem of the type mentioned above would state that under suitable
conditions (\ref{2}), (\ref{3}), (\ref{4}) and (\ref{5)}) are equivalent.

In (\ref{4}) and (\ref{5)}) the functions $u(x)$ and $\tilde{v}(x)$,
respectively, can be considered as utility functions, so that as we move to
the dual versions, that is from (\ref{4}) to (\ref{5)}) , the roles of the
distortions are reversed. This discussion should explain roughly why our
main result, Theorem 1, has four main parts.

\section{The upper $(u_{0},v_{0})$-ordering and a duality theorem}

We start with two functions, $u_{0},v_{0}$ which define our basic stochastic
ordering, following a few definitions. Throughout the paper "increasing"
will mean non-decreasing (unless otherwise stated), and similarly for
"decreasing". All functions will be of bounded variation on compact
intervals and integration is Lebesgue-Stieltjes. Unless otherwise stated,
when we integrate with respect to a measure defined by a right continuous
function the integral will be extended to intervals of the form (a,b], and
when with respect to a left continuous function to intervals [a,b).

\begin{defn}
The pair of functions $(u_{0},v_{0})$ is called \textbf{a standard pair} if

(i) $u_{0}:\boldsymbol{R}\rightarrow \boldsymbol{R}$ is increasing and left
continuous,

(ii) $v_{0}:[0,1]\rightarrow \lbrack 0,1]$ is increasing, right continuous,
and $v_{0}(0)=0,v(1^{-})=1$.\bigskip 
\end{defn}

\begin{defn}
For a function $u_{0}:\boldsymbol{R}\rightarrow \boldsymbol{R}$ the class of 
$u_{0}$-concave functions on $\boldsymbol{R}$ is defined as the class of
functions $\{u:\boldsymbol{R}\rightarrow \boldsymbol{R}\}$ such that 
\begin{equation*}
u(x)=\int_{-\infty }^{x}k(s)du_{0}(s).
\end{equation*}%
for some bounded decreasing function $k(x)$ on $\boldsymbol{R}$. Similarly
define the class of $v_{0}$-concave functions, $\{v:[0,1]\rightarrow \lbrack
0,1]\}$ as those for which 
\begin{equation*}
v(\alpha )=\int_{0}^{\alpha }\tilde{k}(t)dv_{0}(t),
\end{equation*}%
for some bounded decreasing $\tilde{k}$ on $[0,1]$.\bigskip 
\end{defn}

We can interpret Definition 4 as saying that $k(x)$ is the Radon-Nikodym
derivative of $u$ with respect to $u_{0}$ 
\begin{equation*}
k(x)=\frac{du}{du_{0}}
\end{equation*}%
and is decreasing; similarly for $v$ and $v_{0}$.

If we add the condition that $u_{0}$ and $u$ are increasing (see Definition
3), we have that $k(x)\geq 0,x\in \boldsymbol{R}$. We have a similar
interpretation for $\tilde{k}(\alpha )\geq 0,\;\alpha \in \lbrack 0,1]$.

For differentiable $u(x)$ and $u_{0}(x),$ $u_{0}(x)>0$, $u$ being $u_{0}$%
-concave is equivalent to $\dfrac{u^{\prime }}{{u_{0}^{\prime }}}$
decreasing, which in turn, if the functions are twice differentiable and the
second derivatives non-zero, is equivalent to 
\begin{equation}
-\frac{u^{\prime \prime }}{u^{\prime }}\geq -\frac{u_{0}^{\prime \prime }}{%
u_{0}^{\prime }}  \label{6}
\end{equation}%
When $u$ is a utility function, $-\dfrac{u^{\prime \prime }(x)}{u^{\prime
}(x)}$ is the measure of \textit{absolute risk aversion}. By the Arrow-Pratt
Theorem, (Arrow, 1974; Pratt, 1964), $u$ and $u_{0}$ satisfy (\ref{6}) if
and only if 
\begin{equation}
u(x)=\phi (u_{0}(x))\;\;\text{for some concave increasing function}\;\phi .
\label{7}
\end{equation}%
There are similarly versions of (\ref{6}) and (\ref{7}) for the $v_{0}$%
-concave functions of Definition 4. We can also define $u_{0}$- and $v_{0}$%
-convex functions which will be discussed in Section 3.

The stochastic ordering we introduce in this paper is a special example of
the $(U,V)$-ordering discussed above.

\begin{defn}
Given two cdf's $F_{1}$ and $F_{2}$ and a standard pair $(u_{0},v_{0})$,
according to Definition 3 we define the \textbf{upper }$(u_{0},v_{0})$%
\textbf{-stochastic ordering} $F_{1}\prec ^{(u_{0},v_{0})}F_{2}$ as: 
\begin{equation*}
\int_{-\infty }^{\infty }v_{0}(F_{1}(x))du(x)\geq \int_{-\infty }^{\infty
}v_{0}(F_{2}(x))du(x),
\end{equation*}%
for all $u$ in $U^{0}$, the class of $u_{0}$-concave increasing
functions.\bigskip 
\end{defn}

The main result of the paper is the next theorem which is an example of a
duality referred to above. As presaged it says first that $F_{1}\prec
^{(u_{0},v_{0})}F_{2}$ is equivalent to an ordering, involving $F_{1}^{-1}$
and $F_{2}^{-1},$ in which the roles of $u_{0}$ and $v_{0}$ are reversed and
secondly that there are equivalent utility versions, again using $u_{0}$ and 
$v_{0}$ in reverse \textquotedblleft distortion\textquotedblright\ roles.

\begin{thm}
Let $u_{0},v_{0}$ be a standard pair and $U^{0}$ and $V^{0}$ the $u_{0}$%
-concave $v_{0}$-concave increasing classes, respectively. Let $F_{1}$ and $%
F_{2}$ be cdf's which satisfy the following conditions

\bigskip\ (a) $\int_{-\infty }^{\infty }|u_{0}(x)|dv_{0}(F(x))<\infty $

\bigskip\ (b) $\int_{[0,p)}v_{0}(\alpha )du_{0}(F^{-1}(\alpha ))<\infty $
\qquad for all $p<1.$

\vspace{2mm}\noindent Then the following are equivalent:

\bigskip\ (i) $F_{1}\prec ^{(u_{0},v_{0})}F_{2}$

\bigskip\ (ii) $\int_{0}^{1}u_{0}(F_{1}^{-1}(\alpha ))dv(\alpha )\leq
\int_{0}^{1}u_{0}(F_{2}^{-1}(\alpha ))dv(\alpha ),$ for all $v$ in $V^{0}$

\bigskip\ (iii) $\int_{-\infty }^{\infty }u(x)dv_{0}(F_{1}(x))\leq
\int_{-\infty }^{\infty }u(x)dv_{0}(F_{2}(x))$ for all $u$ in $U^{0}$

\bigskip\ (iv) $\int_{0}^{1}v(\alpha )du_{0}(F_{1}^{-1}(\alpha ))\geq
\int_{0}^{1}v(\alpha )du_{0}(F_{2}^{-1}(\alpha ))$ for all $v$ in $V^{0}.$
\end{thm}

There are also four more equivalent statements which are obtained by
transformation of variables in $(i)-(iv)$. For example $(i)^{\ast }$ is
obtained form $(i)$ by a transformation $\alpha =F(x)$ (see (\ref{CV3*}),
Section 9.2 of Appendix 2) and we then obtain a formula with the same
structure as $(iv)$ but with the zero suffix moved from $u_{0}$ to $v$;
similarly for $(ii)^{\ast }$ to $(iv)^{\ast }$.

\bigskip\ $(i)^{\ast }$ $\;\;\int_{0}^{1}v_{0}(\alpha )du(F_{1}^{-1}(\alpha
))\geq \int_{0}^{1}v_{0}(\alpha )du(F_{2}^{-1}(\alpha ))$ for all $u\in
U^{0} $

\bigskip\ $(ii)^{\ast }$ $\;\;\int_{-\infty }^{\infty
}u_{0}(x)dv(F_{1}(x))\leq \int_{-\infty }^{\infty }u_{0}(x)dv(F_{2}(x))$ for
all $v\in V^{0}$

\bigskip\ $(iii)^{\ast }$ $\;\;\int_{0}^{1}u(F_{1}^{-1}(\alpha
))dv_{0}(\alpha )\leq \int_{0}^{1}u(F_{2}^{-1}(\alpha ))dv_{0}(\alpha ),$
for all $u\in U^{0}$

\bigskip\ $(iv)^{\ast }\;\;\int_{-\infty }^{\infty }v(F_{1}(x))du_{0}(x)\geq
\int_{-\infty }^{\infty }v(F_{2}(x))du_{0}(x)$ for all $v\in V^{0}.$

The Proof of Theorem 1 is given in Appendix 1. At the centre of the proof is
the following \textquotedblleft double distortion\textquotedblright\ version
of statement (\ref{1}), also proved in Appendix 1.\bigskip 

\begin{lem}
Let $(u_{0},v_{0})$ be a standard pair and let $F_{1},F_{2}$ be two cdf's on 
$\boldsymbol{R}$ satisfying the condition (a), above, then the following are
equivalent:

\bigskip\ (i) $\int_{-\infty }^{c}v_{0}(F_{1}(x))du_{0}\geq \int_{-\infty
}^{c}v_{0}(F_{2}(x))du_{0}$\qquad for all $c$ in $\boldsymbol{R}$

\bigskip\ (ii) $\int_{0}^{p}u_{0}(F_{1}^{-1}(\alpha ))dv_{0}\leq
\int_{0}^{p}u_{0}(F_{2}^{-1}(\alpha ))dv_{0}$\qquad for all $p$ in $[0,1).$
\end{lem}

\noindent We show that statements $(i)$ and $(ii)$ of Lemma 1 are equivalent
to statements $(i)$ and $(ii)$, respectively, of Theorem 1.

Note that Theorem 1 shows that the upper $(u_{0},v_{0})$-ordering is
equivalent to two inequalities for the \textquotedblleft expected
utilities\textquotedblright\ under the distortion, namely: 
\begin{equation*}
\mbox{E}(u(X))\leq \mbox{E}(u(Y)),
\end{equation*}%
for all $u\in U^{0}$ where $X\sim v_{0}(F_{1})$, $Y\sim v_{0}(F_{2}),$ and
also:%
\begin{equation*}
\mbox{E}(u_{0}(X))\leq \mbox{E}(u_{0}(Y)),
\end{equation*}%
where $X\sim v(F_{1})$, $Y\sim v(F_{2})$ for all $v\in V^{0}.$\bigskip

\section{A convex version}

A review of the stochastic orderings literature points to several results in
which convex increasing functions are used combined with the survivor
function $G(x)=1-F(x)$. The authors pondered as to whether there are two
rather separate theories or whether the duality theory of Section 2 can be
applied without too much additional labour to obtain a convex version. We
believe that indeed the latter is the case and this section develops such a
result. First, we need to define $u_0$-convexity.

\begin{defn}
A function $u:\boldsymbol{R}\rightarrow \boldsymbol{R}$ is said to be $u_{0}$%
-convex if 
\begin{equation*}
u(x)=\int_{-\infty }^{x}m(t)du_{0}(t)
\end{equation*}%
where $m(\cdot )$ is increasing on $\boldsymbol{R}$. Similarly for $v_{0}$%
-convex functions using an increasing $\tilde{m}$ on $[0,1]$.\bigskip 
\end{defn}

We start with a preamble giving transforms which yield a convex version of
Theorem 1.

If $X\sim F(x)$ is a random variable, then the cdf of $-X$ is%
\begin{equation}
F_{-X}(x)=1-F((-x)^{-})=1-F(-x^{+}),  \label{8}
\end{equation}%
and its inverse cdf is 
\begin{equation*}
F_{-X}^{-1}(\alpha )=-F_{X}^{-1}((1-\alpha )^{+})=-F_{X}^{-1}(1-\alpha ^{-}).
\end{equation*}%
Also, for any standard pair $(u_{0},v_{0})$ define 
\begin{equation}
\tilde{u}_{0}(x)=-u_{0}(-x^{-})),\;\;\;\tilde{v}_{0}(\alpha
)=1-v_{0}(1-\alpha ^{+})  \label{9}
\end{equation}%
and note that $(\tilde{u}_{0},\tilde{v}_{0})$ is still a standard pair.

Next, select one of the expressions used in Lemma 1, e.g. 
\begin{equation*}
\int_{-\infty }^{c}v_{0}(F(x))du_{0}(x),
\end{equation*}
and replace $u_0(x)$ by $\tilde{u}_{0}(x)$, $v_{0}(x)$ by $\tilde{v}_{0}(x)$
and $F(x)$ by $F_{-X}(x)$. Making the transformation $z=-x$, we have 
\begin{equation}
\begin{array}{rcl}
\int_{-\infty }^{c}\tilde{v}_{0}(F_{-X}(x))d\tilde{u}_{0}(x) & = & 
-\int_{-\infty }^{c}\tilde{v}_{0}(1-F((-x)^{-}))du_{0}(-x^{-}) \\ 
&  &  \\ 
& = & -\int_{-c}^{\infty }\tilde{v}_{0}(1-F(z^{-}))du_{0}(z) \\ 
&  &  \\ 
& = & -\int_{-c}^{\infty }(1-v_{0}(1-(1-F(z^{-}))^{+})du_{0}(z) \\ 
&  &  \\ 
& = & -\int_{-c}^{\infty }(1-v_{0}(F(z)))du_{0}(z)%
\end{array}
\label{10}
\end{equation}%
These calculations lead naturally to the following definition.

\begin{defn}
For cdf's $F_{1}$ and $F_{2}$ and standard pair $(u_{0},v_{0})$ define the 
\textbf{lower }$(u_{0},v_{0})-$\textbf{\ ordering }$F_{1}\ \prec
_{(u_{0},v_{0})}F_{2}$ by 
\begin{equation*}
-\int_{-\infty }^{\infty }(1-v_{0}(F_{1}(x))du(x)\geq -\int_{-\infty
}^{\infty }(1-v_{0}(F_{2}(x))du(x)
\end{equation*}%
for all increasing $u_{0}$-convex functions $u$.\bigskip 
\end{defn}

The convex version of Theorem 1 is obtained by applying $\tilde{u}_{0},$ $%
\tilde{v}_{0}$, $F_{-X_{1}}$ and $F_{-X_{2}}$ throughout and then converting
back to statements about $u_{0},v_{0},F_{1}$ and $F_{2}$ and using (\ref{8})
and (\ref{9}). It should be added that after this conversion condition $(a)$
from Theorem 1 remains the same, but $(b)$ changes to $(b)^{\prime }$ below.
A compact way of summarizing the analysis is to say that it is a development
of the statement 
\begin{equation}
F_{X_{1}}\ \prec _{(u_{0},v_{0})}F_{X_{2}}\Leftrightarrow F_{-X_{1}}\ \prec
^{(\tilde{u}_{0},\tilde{v}_{0})}F_{-X_{2}}.  \label{11}
\end{equation}%
This could be taken as an equivalent definition.\bigskip 

\begin{thm}
Let $(u_{0},v_{0})$ be a standard pair and $U_{0}$ and $V_{0}$ the $u_{0}$-, 
$v_{0}$-convex increasing classes, respectively. Let $F_{1}$ and $F_{2}$ be
cdf's which satisfy condition $(a)$ of Theorem 1 together with

\vspace{3mm} (b)$^{\prime }$ $\int_{-\infty }^{\infty
}(1-v_{0}(F(x)))du_{0}(x)<\infty $\vspace{3mm}

\noindent Then the following are equivalent:

\vspace{2mm} (i)$^{\prime }$ $F_{1}\prec _{(u_{0},v_{0})}F_{2}$

\vspace{2mm} (ii)$^\prime$ $\int_0^1 u_0(F_1^{-1}(\alpha))dv(\alpha) \leq
\int_0^1 u_0(F_2^{-1}(\alpha))dv(\alpha),$ for all $v$ in $V_0$

\vspace{2mm} (iii)$^\prime$ $\int_{-\infty}^{\infty} u(x) dv_0(F_1(x)) \leq
\int_{-\infty}^{\infty} u(x) dv_0(F_2(x))$ for all $u$ in $U_0$

\vspace{2mm} (iv)$^{\prime }$ $-\int_{0}^{1}(1-v(\alpha
))du_{0}(F_{1}^{-1}(\alpha ))\geq -\int_{0}^{1}(1-v(\alpha
))du_{0}(F_{2}^{-1}(\alpha ))$ for all $v$ in $V_{0}.$\bigskip 
\end{thm}

Important aspects of Theorem 2 are that $u_{0}$-concave and $v_{0}$-concave
are replaced respectively by $u_{0}$-convex and $v_{0}$-convex and that the
utility version for $v(\alpha )$, namely $(iv)$, has as similar structure to 
$v_{0}(\alpha )$ in Definition 7. The new condition, $(b)^{\prime },$
controls the existence of the integrals as $x\rightarrow \infty $ and uses a
distortion generalisation of the survivor function: $1-v_{0}(F(x))$. This
requires that $u(x)$ does not increase too fast as $x\rightarrow \infty $.
This is to be compared to $(b)$ of Theorem 1 which says that $u(x)$ should
not decrease too fast as $x\rightarrow -\infty $.\vspace{3mm}

\section{A combined $(u_{0},v_{0})$-ordering}

Now we combine Theorems 1 and 2 and the upper and lower $(u_{0},v_{0})$%
-orderings and require that $(a),(b)$ and $(b)^{\prime }$, in those
theorems, all hold. First, let us impose no condition on $u_{0}$ except
bounded variation. Then any such $u_{0}$-concave function $u(x)$ can be
represented as 
\begin{equation*}
u(x)=u_{1}(x)+u_{2}(x)
\end{equation*}%
where $u_{1}(x)$ is $u_{0}$-concave increasing and $-u_{2}(x)=u_{3}(x)$ is $%
u_{0}$-convex increasing. This is established by breaking $k(x)$ (see
Definition 4) into non-negative and non-positive parts: $%
k(x)=k^{+}(x)+k^{-}(x)$. Then if we have inequalities involving integrals $%
du_{1}$ and reverse inequalities involving $du_{3}$, \emph{with the same
integrand,} we can achieve bounds with no extra conditions on $u_{0}$.

\begin{defn}
For cdf's $F_{1}$ and $F_{2}$ and a pair of function $(u_{0},v_{0})$ on $%
\boldsymbol{R}$ and $[0,1]$ respectively of bounded variation define the 
\textbf{double }$(u_{0},v_{0})$-\textbf{ordering} as 
\begin{equation*}
\mbox{ $F_{1}\preceq ${\scriptsize (}{\footnotesize u}$_{{\small 0}},${\footnotesize v}$_{{\small 0}}${\scriptsize ) }$F_{2}$
}\Leftrightarrow F_{1}\ \prec ^{(u_{0},v_{0})}F_{2}\;\;\mbox{and}\;\;F_{2}\
\prec _{(u_{0},v_{0})}F_{1}
\end{equation*}%
\bigskip 
\end{defn}

Motivated by the above discussion we have

\begin{lem}
For cdf's $F_{1}$ and $F_{2}$ and a pair of functions $(u_{0},v_{0})$ as in
Definition 8, and satisfying (b) and (b)$^{\prime }$ from Theorems 1 and 2, $%
F_{1}\preceq ${\scriptsize (}{\footnotesize u}$_{{\small 0}},${\footnotesize %
v}$_{{\small 0}}${\scriptsize ) }$F_{2}$ is equivalent to the statement 
\begin{equation}
\int_{-\infty }^{\infty }v_{0}(F_{1}(x))du(x)\leq \int_{-\infty }^{\infty
}v_{0}(F_{2}(x))du(x),  \label{12}
\end{equation}%
for all $u_{0}$-concave functions $u$ (not necessarily increasing).\bigskip 
\end{lem}

\noindent Proof. To establish the inequality in the Lemma we add the
inequalities in definitions of upper and lower orderings which, because of
the reversals, are in the right direction. It is important too that assuming
both $(b)$ and $(b)^{\prime }$ gives the existence of the relevant
integrals. To establish the converse we can make $k^{+}(x)$ and $k^{-}(x)$
in the construction alternatively zero.\vspace{2mm}

Drawing on similar arguments to those for Theorems 1 and 2 we can establish
the following.

\begin{thm}
Let $F_{1}$ and $F_{2}$ be two cdf's and let $(u_{0},v_{0})$ be a standard
pair. Assume $(a),(b)$ and $(b)^{\prime }$ from Theorems 1 and 2 hold; then
the following are equivalent\vspace{2mm}

(i)$^{\prime \prime }$ $F_{1}\ \prec ${\scriptsize (}{\footnotesize u}$_{%
{\small 0}},${\footnotesize v}$_{{\small 0}}${\scriptsize ) }$F_{2}$\vspace{%
2mm}

(ii)$^{\prime \prime }$ $\int_{0}^{1}u_{0}(F_{1}^{-1}(\alpha ))dv(\alpha
)\leq \int_{0}^{1}u_{0}(F_{2}^{-1}(\alpha ))dv(\alpha ),$ for all $v_{0}$%
-concave $v$\vspace{2mm}

(iii)$^{\prime \prime }$ $\int_{-\infty }^{\infty }u(x)dv_{0}(F_{1}(x))\leq
\int_{-\infty }^{\infty }u(x)dv_{0}(F_{2}(x))$ for all $u_{0}$-concave $u$ 
\vspace{2mm}

(iv)$^{\prime \prime }$ $\int_{0}^{1}v(\alpha )du_{0}(F_{1}^{-1}(\alpha
))\leq \int_{0}^{1}v(\alpha )du_{0}(F_{2}^{-1}(\alpha ))$ for all $v_{0}$%
-concave $v$\vspace{2mm}
\end{thm}

\section{Some Examples}

\subsection{Both $u_{0}$ and $v_{0}$ are the identity}

Let $u_{0}(x)=x$ for all $x\in \boldsymbol{R}$ and $v_{0}(\alpha )=\alpha $
for all $\alpha \in $ $[0,1]$. Then $u_{0}$-concave means concave on $%
\boldsymbol{R}$ and\ $v_{0}$-concave means concave on $[0,1]$. Similarly for 
$u_{0}$-convex and $v_{0}$-convex. The upper $(u_{0},v_{0})$-ordering $%
F_{1}\prec ^{(u_{0},v_{0})}F_{2}$ is 
\begin{equation}
\int_{-\infty }^{\infty }F_{1}(x)du(x)\geq \int_{-\infty }^{\infty
}F_{2}(x)du(x)  \label{13}
\end{equation}%
for all increasing concave $u(\cdot )\ $for which the integrals exist.
Condition $(a)$ of Theorem 1 means the existence of the expected values $%
E(X),E(Y)$, where $X\sim $ $F_{1}$ and $Y\sim F_{2}$. In this case some of
the equivalent statements of Theorems 1, 2 and 3 are well-known, but others
are not easily found in the literature.

Theorem 1 states that (\ref{13}) is equivalent to 
\begin{equation}
\int_{-\infty }^{\infty }u(x)dF_{1}(x)\leq \int_{-\infty }^{\infty
}u(x)dF_{2}(x)  \label{14}
\end{equation}%
for all increasing concave $u(\cdot ).$ This is known as the increasing
concave ordering: $F_{1}\leq _{icv}F_{2}$\ (Shaked and Shanthikumar, 2007).
For continuous cdf's the equivalence of (\ref{13}) and (\ref{14}) is
well-known.

Taking $u_{x}(z)=z$ if $z\in (-\infty ,x],$ and $u_{x}(z)=x$ otherwise, for
all $x\in R$, (\ref{13}) is equivalent to 
\begin{equation}
\int_{-\infty }^{x}F_{1}(z)dz\geq \int_{-\infty }^{x}F_{2}(z)dz\qquad 
\mbox
{for all}\;\;x\in R.  \label{15}
\end{equation}%
This is the Second Order Stochastic Dominance ($F_{1}\leq _{SSD}F_{2}$)
ordering. The equivalence of $\leq _{SSD}$ and $\leq _{icv\text{ }}$is also
well-known, see Muller and Stoyan (2002).

Further, by Theorem 1, (\ref{13}) is also equivalent to 
\begin{equation}
\int_{0}^{1}F_{1}^{-1}(\alpha )dv(\alpha )\leq \int_{0}^{1}F_{2}^{-1}(\alpha
)dv(\alpha )  \label{16}
\end{equation}%
for all increasing concave $v(\cdot ),$ and also to%
\begin{equation}
\int_{-\infty }^{\infty }v(F_{1}(x))dx\geq \int_{-\infty }^{\infty
}v(F_{2}(x))dx  \label{17}
\end{equation}%
for all increasing concave $v(\cdot ).$ This equivalence seems to be new.

By taking $v_{p}(\alpha )=\alpha $ if $\alpha \in \lbrack 0,p],$ and $%
v_{p}(\alpha )=p$ otherwise, for all $\alpha \in \lbrack 0,1]$, (\ref{16})
becomes 
\begin{equation}
\int_{0}^{p}F_{1}^{-1}(\alpha )d\alpha \leq \int_{0}^{p}F_{2}^{-1}(\alpha
)d\alpha \qquad \;\;\mbox{for
all}\;\;p\in \lbrack 0,1]  \label{18}
\end{equation}%
which is a generalization of the Lorenz ordering. The equivalence of (\ref%
{18}) and $\leq _{icv}$ can be found in several papers, at least for special
cases. Recently Sordo and Ramos (2007) have proved it under general
conditions.\bigskip

With $u_{0}$ and $v_{0}$ the identity, the lower $(u_{0},v_{0})$-ordering $%
F_{1}\leq _{(u_{0},v_{0})}F_{2}$ can be expressed as: 
\begin{equation}
\int_{-\infty }^{\infty }G_{1}(x)du(x)\leq \int_{-\infty }^{\infty
}G_{2}(x)du(x)  \label{19}
\end{equation}%
for all increasing convex $u(\cdot )$, where $G_{j}(x)=1-F_{j}(x),$ $%
(j=1,2). $ By Theorem 2, (\ref{19}) is equivalent to

\begin{equation}
\int_{-\infty }^{\infty }u(x)dF_{1}(x)\leq \int_{-\infty }^{\infty
}u(x)dF_{2}(x)  \label{20}
\end{equation}%
for all increasing convex $u(\cdot ).$ This is known as the increasing
convex ($\leq _{icx}$) ordering. It is also equivalent to 
\begin{equation}
\int_{0}^{1}F_{1}^{-1}(\alpha )dv(\alpha )\leq \int_{0}^{1}F_{2}^{-1}(\alpha
)dv(\alpha )  \label{21}
\end{equation}%
for all increasing convex $v(\cdot ).$ The equivalence of $\leq _{icx}$(\ref%
{20}) and (\ref{21}) is a well-known result, see Shaked and Shanthikumar
(2007), Theorem 4.A.4.

By taking $v_{p}(\alpha )=\alpha $ if $\alpha \in (p,1],$ $v_{p}(\alpha )=p$
otherwise, for all $\alpha \in \lbrack 0,1]$, (\ref{21}) becomes also
equivalent to 
\begin{equation}
\int_{p}^{1}F_{1}^{-1}(\alpha )d\alpha \leq \int_{p}^{1}F_{2}^{-1}(\alpha
)d\alpha \qquad \;\;\mbox{for all}\;\;p\in \lbrack 0,1]  \label{22}
\end{equation}%
and the equivalence of $\leq _{icx}$ and (\ref{22}) is Theorem 4.A.3 of
Shaked and Shanthikumar (2007).

\subsection{Majorization}

Majorization\ is an ordering $\prec $ of real vectors which satisfies the
Dalton-Pigou Principle of Transfers of wealth and was brought to the fore by
Marshall and Olkin's (1979) fundamental book. Given the vectors $\boldsymbol{%
x}=(x_{1},x_{2},...,x_{n})$ and $\boldsymbol{y}=(y_{1},y_{2},...,y_{n})$ let 
$\ x_{(1)}\geq x_{(2)}\geq ...\geq x_{(n)}$ and\ $y_{(1)}\geq y_{(2)}\geq
...\geq y_{(n)}$ be the rearranged coordinates; one of several equivalent
definitions of $\boldsymbol{y}\prec \boldsymbol{x}$ is%
\begin{eqnarray*}
\sum_{i=1}^{k}x_{(i)} &\geq &\sum_{i=1}^{k}y_{(i)}\qquad k=1,...,n-1 \\
\sum_{i=1}^{n}x_{(i)} &=&\sum_{i=1}^{n}y_{(i)}
\end{eqnarray*}

\noindent Removing the \textquotedblleft equal means\textquotedblright\
condition (bottom line) defines two extensions of this ordering, that are
perhaps less well known: \textit{lower weak majorization} $\boldsymbol{y}%
\prec _{w}\boldsymbol{x:}$%
\begin{equation*}
\sum_{i=1}^{k}x_{(i)}\geq \sum_{i=1}^{k}y_{(i)}\qquad k=1,...,n\ 
\end{equation*}%
and \textit{upper weak majorization }$\boldsymbol{y}\prec ^{w}\boldsymbol{x}$
\begin{equation}
\sum_{i=n-k}^{n}x_{(i)}\leq \sum_{i=n-k}^{n}y_{(i)}\qquad k=0,1,...,n-1
\label{(uwm)}
\end{equation}

\noindent Our theory gives results well known or easily derived directly.
Let $X$ and $Y$ be random variables on finite supports on the line, such
that $\Pr (X=x_{i})=$ $\Pr (Y=y_{i})=1/n$ for all $i$. The upper $%
(u_{0},v_{0})$-ordering $F_{1}\prec ^{(u_{0},v_{0})}F_{2}$ with $u_{0}$ the
identity and $v_{0}$ the identity is upper weak majorization \textit{\ }$%
\boldsymbol{y}\prec ^{w}\boldsymbol{x}$ (in the reverse ordering of the
vectors). By Lemma 1, $\prec ^{(u_{0},v_{0})}$ is defined by $%
\int_{0}^{p}F_{1}^{-1}(\alpha ))d\alpha \leq \int_{0}^{p}F_{2}^{-1}(\alpha
)d\alpha $ for all $0<p\leq 1,$ which is precisely (\ref{(uwm)}). An
increasing concave function $u(x)$ yields positive decreasing increments $%
u(x_{(n-i+1)})-u(x_{(n-i)})$ and $u(y_{(n-i+1)})-u(y_{(n-i)})$ for $%
i=1,...,n-1.$ Similarly $v$ increasing concave means that the increments $%
b_{(i)}=v(\frac{i}{n})-v(\frac{i-1}{n})$ yield a positive decreasing
sequence,\textit{\ }and the equivalent statements (i) to (iv) in Theorem 1
can be translated into equivalent statements (\ref{24}) to (\ref{(27)})
below. Equivalence with (\ref{24}) to (\ref{26}) can be found in Marshall
and Olkin (1979) whereas (\ref{(27)}) can be easily obtained from the
definition

\noindent (i)%
\begin{equation}
\sum_{i=1}^{n-1}\frac{i}{n}(u(x_{(n-i)})-u(x_{(n-i+1)})\geq \sum_{i=1}^{n-1}%
\frac{i}{n}(u(y_{(n-i)})-u(y_{(n-i+1)}))  \label{24}
\end{equation}%
for any increasing concave real function $u(\cdot );$\vspace{2mm}

\noindent (ii)%
\begin{equation}
\sum_{i=1}^{n}b_{\left( i\right) }x_{(i)}\leq \sum_{i=1}^{n}b_{\left(
i\right) }y_{(i)}  \label{25}
\end{equation}%
for any decreasing sequence $1\geq b_{\left( 1\right) }\geq b_{\left(
2\right) }\geq ...\geq b_{\left( n\right) }\geq 0$;\vspace{2mm}

\noindent (iii)%
\begin{equation}
\sum_{i=1}^{n}u(x_{i})\leq \sum_{i=1}^{n}u(y_{i})  \label{26}
\end{equation}%
for all increasing concave real functions $u(\cdot )$;\vspace{2mm}

\noindent (iv)%
\begin{equation}
\sum_{i=1}^{n}v(\frac{i}{n})(x_{(n-i)}-x_{(n-i+1)})\geq \sum_{i=1}^{n}v(%
\frac{i}{n})(y_{(n-i)}-y_{(n-i+1)})  \label{(27)}
\end{equation}%
for all increasing concave real functions $v(\cdot ),$ where $%
x_{(0)}=y_{(0)}=K$ is any real number.\vspace{2mm}

Similarly, for uniform distributions with finite supports on the real line,
the lower ($u_{0},v_{0})$-ordering $F_{1}\prec _{(u_{0},v_{0})}F_{2}$ when $%
u_{0}$ and $v_{0}$ are the identity becomes $\boldsymbol{x}\prec _{w}%
\boldsymbol{y}$ and the equivalence of $(i),(ii),(iii)$ and $(iv)$ of
Theorem 2 gives equivalent statements in the theory of lower weak
majorization.

It is well known that standard majorization $\boldsymbol{x}\prec \boldsymbol{%
y}$ holds iff $\boldsymbol{x}\prec ^{w}\boldsymbol{y}$ and $\boldsymbol{x}%
\prec _{w}\boldsymbol{y}$. Thus $\boldsymbol{x}\prec \boldsymbol{y}$ is a
special case of the double ordering of Theorem 3 with $u_{0},v_{0}$ both the
identity.

It may be interesting to consider extensions when $u_{0}$ is a different
increasing function. For instance, $u_{0}(x)=\log x$ yields the ordering
known as log-majorization (including the weak versions).\ Theorems 1, 2 and
3 in this case provide non-trivial new results.\bigskip

\subsection{Social welfare functionals}

The risk functional $\int_{-\infty }^{\infty }u_{0}(x)dF(x)$ is the main
mathematical entity of Expected Utility (EU) theory and the theory of
decision-making under risk. There is a close formal relationship between
that literature and the literature on income distribution. From a
mathematical point of view it is clear there is an extensive duality theory
between the "forward" theory of utility, and the reverse theory, the theory
of welfare. We mention a few implications of our results for these theories.
On the side of the welfare literature are the so-called "rank-dependent
social welfare functionals" and the related rank-dependent expected utility
theory (RDEU). This goes back to Quiggin (1982) who developed the RDEU model
based on 
\begin{equation}
\int_{-\infty }^{\infty }u_{0}(x)df_{0}(F(x))  \label{(28)}
\end{equation}%
where $f_{0}:[0,1]\rightarrow \lbrack 0,1]$ is a strictly increasing and
continuous distortion function, called a \textquotedblleft
perception\textquotedblright , for which $f_{0}(0)=0$ and $f_{0}(1)=1$ and $%
u_{0}$ is a strictly increasing utility. When $f_{0}$ is the identity
function the RDEU theory reduces to EU. With the Gini index as a prominent
exception, the most common inequality measures can be interpreted in a
social welfare framework formally equivalent to the EU model. When $u_{0}$
is the identity, the RDEU model yields the Yaari welfare function 
\begin{equation}
W(F)=\int_{-\infty }^{\infty }xdf_{0}(F(x)).  \label{(29)}
\end{equation}%
The most popular social welfare function of the form $W(F)$ is the \textit{%
S-Gini function}, where $f_{0}(p)=p^{\rho }$ with $\rho >1$ (Donaldson and
Weymark, 1980; Yitzhaki, 1983). Note that the classical Gini index is
associated to the S-Gini social welfare function with $\rho =2$. The
parameter $\rho $ can be seen as a measure of inequality aversion. Note that
(\ref{(29)}) can be rewritten as

\begin{equation}
W(F)=-\int_{-\infty }^{\infty }xdv_{0}(\bar{F}(x)),  \label{(30)}
\end{equation}%
where $\bar{F}(x)=1-F(x)$ and $v_{0}(\alpha )=1-f_{0}(1-\alpha )$, so that $%
0\leq v_{0}(\alpha )\leq 1,$ and $v_{0}$ is increasing. Increasing $v_{0}$%
-concave functions are all the increasing $f_{0}$-convex ones. A simple
change of variables yields a different expression for Yaari's functional: 
\begin{equation}
W(F)=-\int_{-\infty }^{\infty }xdv_{0}(\bar{F}(x))=\int_{0}^{1}F^{-1}(\alpha
)df_{0}(\alpha ).  \label{(31)}
\end{equation}%
We now state two simple corollaries of Theorem 1.

\begin{cor}
\textbf{\qquad }If $u_{0}$ is the identity and $f_{0}$ an increasing
function, under the assumptions of Theorem 1 the following two statements
are equivalent: 
\begin{equation}
\int_{-\infty }^{\infty }u(x)df_{0}(F_{1}(x))\leq \int_{-\infty }^{\infty
}u(x)df_{0}(F_{2}(x)),\;\;\mbox{for all}\;u\text{ increasing and concave} 
\tag{S$_{\text{1}}$}
\end{equation}%
\begin{equation}
\int_{0}^{1}F_{1}^{-1}(\alpha ))df(\alpha )\leq
\int_{0}^{1}F_{2}^{-1}(\alpha )df(\alpha ),\;\text{for all}\;\text{%
increasing }f\text{ that are }f_{0}\text{-convex.}  \tag{S$_{\text{2}}$}
\end{equation}%
\vspace{2mm}
\end{cor}

The first of these statements say that $F_{2}$ is preferred to $F_{1}$, in
the expected utility sense, by all risk-averse decision-makers with
\textquotedblleft perception\textquotedblright\ defined by $f_{0}$. The
second says that $F_{2}$ is preferred to $F_{1}$, in the Yaari sense, by
decision-makers whose inequality aversion is greater than $f_{0}$.

The equivalence of (S$_{\text{1}}$) and (S$_{\text{2}}$) can be extended, as
a consequence of the following corollary.\vspace{2mm}

\begin{cor}
\textbf{\qquad }Take $u_{0}$ and $f_{0}$ to be increasing functions. When
the assumptions of Theorem 1 are satisfied, the following two statements are
equivalent: 
\begin{equation}
\int_{-\infty }^{\infty }u(x)df_{0}(F_{1}(x))\leq \int_{-\infty }^{\infty
}u(x)df_{0}(F_{2}(x)),\;\;\mbox{for all}\;u\ \text{ increasing\ and\ }u_{0}%
\text{-concave,}  \tag{S$_{\text{3}}$}
\end{equation}%
\begin{equation}
\int_{-\infty }^{\infty }u_{0}(x)df(F_{1}(x))\leq \int_{-\infty }^{\infty
}u_{0}(x)df(F_{2}(x)),\;\;\mbox{for all increasing}\;\;f\text{ \ that are }%
f_{0}\text{-convex.}  \tag{S$_{\text{4}}$}
\end{equation}%
\vspace{2mm}
\end{cor}

Inequality (S$_{\text{3}}$) expresses the preference of $F_{2}$ over $F_{1}$
in a utility sense, for all decision-makers with a given perception\ $f_{0}$
whose risk aversion is greater than $u_{0}$; inequality (S$_{\text{4}}$)
states the preference of $F_{2}$ over $F_{1}$ in a RDEU sense by all
decision-makers with utility $u_{0}$ and greater inequality aversion than $%
f_{0}$. To the best of these authors' knowledge, the equivalence of (S$_{%
\text{3}}$) and (S$_{\text{4}}$) does not appear in the literature.

\section{Discussion}

Stochastic orderings are an attractive way of summarizing preferences
between distributions such as in comparing portfolios, assessing risk in
insurance, in individual decision-making and in the study of income
distribution and welfare. Our starting point in Definition 1 is to stress
the use of partial orderings, that is group preferences, via the $(U,V)$%
-formulation.

Our version of this integral stochastic ordering, in Definition 5, captures
the preference not simply of a single subject but those of a group of
subjects each with a private utility and each at least as risk averse as a
base subject represented by $u_{0}$. The duality theory says that this group
defines a dual group described by utilities attached to the quantile
function, with its own base utility $v_{0}$. Members of this dual group are
at least as risk averse as the base subject in the dual realm. Moreover, the
utility function for the base subject in the first group provides a
(probability) distortion in the dual theory and vice versa.\bigskip

\section{Appendix 1}

\subsection{Proof of Theorem 1}

The proof is in two stages. The first is to show that Lemma 1 holds at
\textquotedblleft crossing points\textquotedblright\ of $F_{1}$ and $F_{2}$,
the second is to extrapolate the result between crossing points.

\begin{defn}
For two cdf's $F_{1}$ and $F_{2}$ a \textbf{crossing interval} is a set $%
[a,b]\subset \boldsymbol{R}$ such that there exists $\epsilon >0$ such that
for all $0<\epsilon _{1},\epsilon _{2}<\epsilon $

(i) $F_{1}(a-\epsilon _{1})<F_{2}(a-\epsilon _{1})$

(ii) $F_{1}(b+\epsilon _{2})>F_{2}(b+\epsilon _{2})$

(iii) $\mbox{when}\;\;a<b,\;F_{1}(x)=F_{2}(x),\;$for all $x\in (a,b)$

\noindent In this case we say that \textquotedblleft $F_{1}$ up-crosses $%
F_{2}$\textquotedblright . If the roles of $F_{1}$ and $F_{2}$ are reversed
we say that \textquotedblleft $F_{1}$ down-crosses $F_{2}$\textquotedblright
.\vspace{2mm}
\end{defn}

Thus an up-crossing (down-crossing) \textit{point} $x_{0}$ is such that $%
F_{1}(x_{0}^{-})\leq F_{2}(x_{0}^{-})\leq F_{2}(x_{0})\leq F_{1}(x_{0})$ ($%
F_{2}(x_{0}^{-})\leq F_{1}(x_{0}^{-})\leq F_{1}(x_{0})\leq F_{2}(x_{0}))$.
We can similarly define crossing intervals $[\alpha ,\beta ]\subset \lbrack
0,1]$ for $F_{1}^{-1}$ and $F_{2}^{-1}$. With care, we can make the crossing
intervals match up: if $x_{0}$ is an up-crossing point of $(F_{1},F_{2})$,
then $[F_{2}(x_{0}^{-}),F_{2}(x_{0})]$ is a down-crossing interval for $%
(F_{1}^{-1},F_{2}^{-1})$ and similarly for the converse.

\begin{lem}
If $x_{0}$ is a crossing point for the pair $(F_{1},F_{2}),$ then%
\begin{equation*}
\int_{-\infty }^{x_{0}}v_{0}(F_{1}(x))du_{0}(x)\geq \int_{-\infty
}^{x_{0}}v_{0}(F_{2}(x))du_{0}(x)
\end{equation*}%
implies%
\begin{equation*}
\int_{0}^{p}u_{0}(F_{1}^{-1}(\alpha ))dv_{0}(\alpha )\leq
\int_{0}^{p}u_{0}(F_{2}^{-1}(\alpha ))dv_{0}(\alpha )
\end{equation*}%
for all $p\in \lbrack F_{2}(x_{0}^{-}),F_{2}(x_{0})]$ if $F_{1}$ up-crosses $%
F_{2}$ and for all $p\in \lbrack F_{1}(x_{0}^{-}),F_{1}(x_{0})]$ if $F_{1}$
down-crosses $F_{2}$. Similarly, given an up-crossing value $\alpha _{0}$ of 
$(F_{1}^{-1},F_{2}^{-1})$%
\begin{equation*}
\int_{0}^{\alpha _{0}}u_{0}(F_{1}^{-1}(\alpha ))dv_{0}(\alpha )\leq
\int_{0}^{\alpha _{0}}u_{0}(F_{2}^{-1}(\alpha ))dv_{0}(\alpha )
\end{equation*}%
implies%
\begin{equation*}
\int_{-\infty }^{c}v_{0}(F_{1}(x))du_{0}(x)\geq \int_{-\infty
}^{c}v_{0}(F_{2}(x))du_{0}(x)
\end{equation*}%
for all $c$ in $[F_{2}^{-1}(\alpha _{0}),F_{2}^{-1}(\alpha _{0}^{+})]$. The
same for a down-crossing value, changing the statements appropriately.
\end{lem}

\noindent Proof. Change of variables is the key to the proof. Thus, by Lemma
4 in Appendix 2, and the discussion there, applied first to $F_{1}(x)$ and
then to$\ F_{2}(x)$, we have 
\begin{gather}
\int_{-\infty
}^{x_{0}}v_{0}(F_{1}(x))du_{0}(x)+\int_{0}^{p}u_{0}(F_{1}^{-1}(\alpha
))dv_{0}(\alpha )=  \label{(CV1)} \\
=v_{0}(F_{1}(x_{0}))u_{0}(x_{0})+v_{0}(p)u_{0}(F_{1}^{-1}(p))-v_{0}(F_{1}(x_{0}))u_{0}(F_{1}^{-1}(p))
\notag
\end{gather}%
and%
\begin{gather}
\int_{-\infty
}^{x_{0}}v_{0}(F_{2}(x))du_{0}(x)+\int_{0}^{p}u_{0}(F_{2}^{-1}(\alpha
))dv_{0}(\alpha )=  \label{CV2} \\
=v_{0}(F_{2}(x_{0}))u_{0}(x_{0})+v_{0}(p)u_{0}(F_{2}^{-1}(p))-v_{0}(F_{2}(x_{0}))u_{0}(F_{2}^{-1}(p))
\notag
\end{gather}

Let $x_{0}$ be an up-crossing point for the pair $(F_{1},F_{2})$ and let $p$ 
$\in $ $[F_{2}(x_{0}^{-}),F_{2}(x_{0})].$ For the purpose of proving the
first implication in Lemma 3, we can assume, without loss of generality,
that any crossing interval reduces to just a point, since open intervals on
which $F_{1}(x)=F_{2}(x)$ (part (iii) of Definition 9) only contribute by
adding a constant to the left hand sides of the last identities. By the same
argument, for all $F_{1}(x_{0}^{-})\leq F_{2}(x_{0}^{-})\leq \alpha \leq
F_{2}(x_{0})\leq F_{1}(x_{0})$ we have $F_{1}^{-1}(\alpha
)=F_{2}^{-1}(\alpha )=x_{0},$ so without loss of generality we can assume,
similarly, that $x_{0}=F_{1}^{-1}(p)=F_{2}^{-1}(p).${\LARGE \ }Then the
right hand sides of (\ref{(CV1)}) and (\ref{CV2}) become equal and the first
implication in Lemma 3 is true for up-crossing. The proof for down-crossing
points is similar. Furthermore, the proof of the converse is now
straightforward.\bigskip 

The value of Lemma 3 is to highlight \textquotedblleft
good\textquotedblright\ points where it is straightforward to prove the
equivalence in Lemma 1. It is a little more straightforward to prove the
reverse implication (ii) $\Rightarrow $ (i) in Lemma 1 first. Thus, assume
that (ii) in Lemma 1 holds for all $p$. Then for any $x_{0}$ which belongs
to a crossing interval the inequality (i) in Lemma 1 holds, as just shown.
We now need to extend the proof essentially to the regions between crossing
intervals.

Thus, suppose that for a given $x_{0}$ there is no $p$ such that $(x_{0},p)$
is a crossing pair. Let\ 
\begin{equation*}
x_{1}=\sup \{x^{\prime }:\;\mbox{such that}\;x^{\prime }<x_{0}\;\text{and}%
\;x^{\prime \text{ }}\text{is a\ crossing point}\}
\end{equation*}%
\medskip then $x_{1}$ belongs to a crossing interval and assume 
\begin{equation*}
\begin{array}{rcl}
\int_{-\infty }^{x_{1}}v_{0}(F_{1}(x))du_{0} & \geq & \int_{-\infty
}^{x_{1}}v_{0}(F_{2}(x))du_{0}%
\end{array}%
\end{equation*}%
from the assumption. If at this crossing interval $F_{1}$ up-crosses $F_{2}$
then 
\begin{equation*}
\begin{array}{r}
v_{0}(F_{1}(x))\geq v_{0}(F_{1}(x))%
\end{array}%
\end{equation*}%
for all $x_{1}\leq x\leq x_{0}$ and 
\begin{equation*}
\begin{array}{lll}
\int_{-\infty }^{x_{0}}v_{0}(F_{1}(x))du_{0} & = & \int_{-\infty
}^{x_{1}}v_{0}(F_{1}(x))du_{0}+\int_{x_{1}}^{x_{0}}v_{0}(F_{1}(x))du_{0} \\ 
&  &  \\ 
& \geq & \int_{-\infty
}^{x_{1}}v_{0}(F_{2}(x))du_{0}(x)+\int_{x_{1}}^{x_{0}}v_{0}(F_{2}(x))du_{0}
\\ 
&  &  \\ 
& = & \int_{-\infty }^{x_{0}}v_{0}(F_{2}(x))du_{0}%
\end{array}%
\end{equation*}%
\medskip Note that possibly $x_{1}=-\infty $, in which case $%
v_{0}(F_{1}(x))\geq v_{0}(F_{2}(x))$ for all $x\leq x_{0}$ and the assertion
remains true.

If $F_{1}$ down-crosses $F_{2}$ at $x_{1},$ let%
\begin{equation*}
x_{2}=\inf \{x^{\prime }:\mbox{such that}\;x^{\prime }>x_{0}\;\text{and}%
\;x^{\prime }\text{ is a\ crossing point}\}.
\end{equation*}%
\medskip so that $x_{2}$ is an up-crossing point for $(F_{1},F_{2}).$ Note
that in this case possibly $x_{2}=+\infty $.

Now assume 
\begin{equation*}
\begin{array}{rcl}
\int_{-\infty }^{x_{2}}v_{0}(F_{1}(x))du_{0} & \geq & \int_{-\infty
}^{x_{2}}v_{0}(F_{2}(x))du_{0}.%
\end{array}%
\end{equation*}%
\medskip Also, since $v_{0}(F_{1}(x))\leq v_{0}(F_{2}(x))$ for all $%
x_{0}\leq x\leq x_{2}$, we obtain 
\begin{equation*}
\begin{array}{lll}
\int_{-\infty }^{x_{0}}v_{0}(F_{1}(x))du_{0} & = & \int_{-\infty
}^{x_{2}}v_{0}(F_{1}(x))du_{0}-\int_{x_{0}}^{x_{2}}v_{0}(F_{1}(x))du_{0} \\ 
&  &  \\ 
& \geq & \int_{-\infty
}^{x_{2}}v_{0}(F_{2}(x))du_{0}-\int_{x_{0}}^{x_{2}}v_{0}(F_{2}(x))du_{0} \\ 
&  &  \\ 
& = & \int_{-\infty }^{x_{0}}v_{0}(F_{2}(x))du_{0}%
\end{array}%
\end{equation*}%
The forward implication, (i) $\Rightarrow $ (ii) in Lemma 1 follows on
similar lines, starting with an arbitrary $p\in \lbrack 0,1]$. This
completes the proof of Lemma 1.\bigskip

The next step in the proof of Theorem 1 involves a mixing argument. We start
with (i) in Lemma 1 and claim that for any given $c\in \boldsymbol{R}$%
\medskip\ 
\begin{equation*}
\begin{array}{rcl}
\int_{-\infty }^{c}v_{0}(F_{1}(x))du_{0}(x) & \geq & \int_{-\infty
}^{c}v_{0}(F_{2}(x))du_{0}(x) \\ 
& \Leftrightarrow &  \\ 
\int_{-\infty }^{\infty }\int_{-\infty }^{c}v_{0}(F_{1}(x))du_{0}(x)d\mu (c)
& \geq & \int_{-\infty }^{\infty }\int_{-\infty
}^{c}v_{0}(F_{2}(x))du_{0}(x)d\mu (c),%
\end{array}%
\end{equation*}%
for all non-negative bounded $\sigma $-finite measures $d\mu (c)$ on $%
\boldsymbol{R}$. Next, introducing the indicator function $\mathbb{I}%
_{(-\infty ,c)}(x),$ reversing the integrals and using Fubini's Theorem,
which holds because of the boundedness of $d\mu (c)$ and condition $(a)$, we
have 
\begin{equation*}
\begin{array}{rl}
\int_{-\infty }^{\infty }\left[ \int_{-\infty }^{\infty }v_{0}(F_{1}(x)){%
\mathbb{I}}_{(-\infty ,c)}(x)du_{0}(x)\right] d\mu (c) & \geq \\ 
& \int_{-\infty }^{\infty }\left[ \int_{-\infty }^{\infty }v_{0}(F_{2}(x)){%
\mathbb{I}}_{(-\infty ,c)}(x)du_{0}(x)\right] d\mu (c) \\ 
& \Leftrightarrow \\ 
\int_{-\infty }^{\infty }v_{0}(F_{1}(x))\left[ \int_{-\infty }^{\infty }{%
\mathbb{I}}_{(-\infty ,c)}(x)d\mu (c)\right] du_{0}(x) & \geq \\ 
& \int_{-\infty }^{\infty }v_{0}(F_{2}(x))\left[ \int_{-\infty }^{\infty }{%
\mathbb{I}}_{(-\infty ,c)}(x)d\mu (c)\right] du_{0}(x) \\ 
& \Leftrightarrow \\ 
\int_{-\infty }^{\infty }v_{0}(F_{1}(x))k(x)du_{0}(x)\geq & \int_{-\infty
}^{\infty }v_{0}(F_{2}(x))k(x)du_{0}(x) \\ 
& \Leftrightarrow \\ 
\int_{-\infty }^{\infty }v_{0}(F_{1}(x))du(x)\geq & \int_{-\infty }^{\infty
}v_{0}(F_{2}(x))du(x),%
\end{array}%
\end{equation*}%
where $k(x)=\int_{-\infty }^{\infty }{\mathbb{I}}_{(-\infty ,c)}(x)d\mu (c)$
is a non-negative decreasing bounded function and we define $u(x)$ so that $%
k(x)du_{0}=du$. But such a $u$ is precisely a $u_{0}$-concave increasing
function satisfying Definition 4. A similar argument applies to statement $%
(ii)$ in Lemma 1 and we obtain the equivalent version. Thus, we have shown
that $(i)\Leftrightarrow (ii)$ in Theorem 1. Note that we use condition $(a)$
to obtain Fubini in this case and a bounded decreasing function $\tilde{k}%
(\alpha ),\;\alpha \in \lbrack 0,1]$, as in Definition 4.

Finally, that condition $(iii)$ is equivalent to $(i)$, and $(ii)$ is
equivalent to $(iv)$ follows from the version of integration by parts in
Section 6.1 of Appendix 2, the discussion therein and conditions $(a)$ and $%
(b)$ in Theorem 1. This is to obtain bounded integrals. This ends the proof
of Theorem 1.\bigskip

\section{Appendix 2}

\subsection{Integration by parts}

Because we use a nonstandard version of the \textquotedblleft integration by
parts\textquotedblright\ theorem, we include a full proof here.

\begin{thm}
\qquad Let U and V denote two real functions of finite variation on each
compact interval of the real line, with U\ left continuous and V right
continuous. Then for each pair of real numbers $a<b$%
\begin{equation*}
\int_{(a,b]}U(x)dV+\int_{[a,b)}V(x)dU=U(b)V(b)-V(a)U(a)
\end{equation*}
\end{thm}

\noindent Proof. Define the measures%
\begin{eqnarray*}
\mu \left\{ \lbrack a,b)\right\} &=&U(b)-U(a) \\
\nu \left\{ (a,b]\right\} &=&V(b)-V(a)
\end{eqnarray*}%
for each pair of real numbers $a<b$. The statement of the theorem is
equivalent to the following

\begin{center}
$%
\begin{array}{l}
\int_{(a,b]}\left( U(x)-U(a)\right) d\nu +\int_{[a,b)}\left(
V(x)-V(a)\right) d\mu = \\ 
\\ 
=U(b)V(b)-V(a)U(a)-V(a)\left( U(b)-U(a)\right) -U(a)\left( V(b)-V(a)\right) ,%
\end{array}%
$
\end{center}

\noindent which is also equivalent to{\normalsize 
\begin{equation*}
\int_{(a,b]}\left( U(x)-U(a)\right) d\nu +\int_{[a,b)}\left(
V(x)-V(a)\right) d\mu =\left( U(b)-U(a)\right) \left( V(b)-V(a)\right) .
\end{equation*}%
}

\noindent We observe that 
\begin{eqnarray*}
\int_{\lbrack a,b)}\left( V(x)-V(a)\right) d\mu &=&\int_{[a,b)}\nu \{y\in
(a,b]:y\leq x\}\mu (dx) \\
&=&\mu \otimes \nu \{(x,y)\in \lbrack a,b)\times (a,b]:y\leq x\}
\end{eqnarray*}%
by Fubini's theorem. Similarly 
\begin{eqnarray*}
\int_{(a,b]}\left( U(x)-U(a)\right) d\nu &=&\int_{(a,b]}\mu \{x\in \lbrack
a,b):x<y\}\nu (dy) \\
&=&\mu \otimes \nu \{(x,y)\in \lbrack a,b)\times (a,b]:x<y\}.
\end{eqnarray*}%
Adding up these two identities, we obtain 
\begin{eqnarray}
\int_{(a,b]}\left( U(x)-U(a)\right) d\nu +\int_{[a,b)}\left(
V(x)-V(a)\right) d\mu &=&\mu \otimes \nu \{[a,b)\times (a,b]\}  \label{29} \\
&=&\left( U(b)-U(a)\right) \left( V(b)-V(a)\right) ,  \notag
\end{eqnarray}%
an equivalent statement to the assert of this theorem.\bigskip

We first apply Theorem 4 taking ($u_{0},v_{0}$) to be a standard pair, and
letting $U(x)=u_{0}(x)$ and $V(x)=v_{0}(F(x))$. Then under condition (a) of
Theorem 1 we obtain 
\begin{eqnarray}
\int_{(-\infty ,c]}u_{0}(x)dv_{0}(F(x))+\int_{(-\infty
,c)}v_{0}(F(x))du_{0}(x) &=&u_{0}(c)v_{0}(F(c))-\lim_{a\rightarrow -\infty
}u_{0}(a)v_{0}(F(a))  \notag \\
&=&u_{0}(c)v_{0}(F(c))  \label{30}
\end{eqnarray}%
Then we apply Theorem 4 again taking $U(\alpha )=u_{0}(F^{-1}(\alpha ))$ and 
$V(\alpha )=v_{0}(\alpha )$ and, again under condition $(a)$ of Theorem 1,
obtain 
\begin{align}
\int_{(0,p]}u_{0}(F^{-1}(\alpha ))dv_{0}(\alpha )+\int_{[0,p)}v_{0}(\alpha
)du_{0}(F^{-1}(\alpha ))& =u_{0}(F^{-1}(p))v_{0}(p)-\lim_{\alpha \rightarrow
0}u_{0}(F^{-1}(\alpha ))v_{0}(\alpha )  \label{31} \\
& =u_{0}(F^{-1}(p))v_{0}(p).  \notag
\end{align}%
\bigskip

\subsection{Change of variables}

A classical mathematical result states.

\begin{thm}
(Change of variables). \ Let $(\Omega ,\Xi ,\mu )$ be a measure space and $%
\varphi :\Omega \rightarrow \boldsymbol{R}$ $\ $a measurable function. For
any Borel set $A$ consider the measure $\mu \varphi ^{-1}(A)=\mu (\varphi
^{-1}(A)).$ Let $f$ be measurable real function on the real line $\mathbf{R}$%
. Then for all $A$%
\begin{equation}
\int_{\varphi ^{-1}(A)}f(\varphi (\omega ))\mu (d\omega )=\int_{A}f(x)\mu
\varphi ^{-1}(dx).  \label{32}
\end{equation}
\end{thm}

We apply this theorem to the sets $A=(-\infty ,a]$ and the function $\varphi
(\alpha )=F^{-1}(\alpha )$ where $F$ is a cdf so that $\Omega =$ $(0,1).$
Recall%
\begin{equation}
F(x)\geq \alpha \iff F^{-1}(\alpha )\leq x\qquad \forall x\in \boldsymbol{R}%
,\forall \alpha \in \lbrack 0,1]  \label{33}
\end{equation}%
hence%
\begin{equation}
F(x)<\alpha \iff F^{-1}(\alpha )>x\qquad \forall x\in \boldsymbol{R},\forall
\alpha \in \lbrack 0,1].  \label{34}
\end{equation}%
It is easy to see that $\varphi ^{-1}(a,b]=(F(a),F(b)]$ for all $a,b\in 
\boldsymbol{R}$. Take $f=u_{0}$\ and $\mu $ the measure on $[0,1]$ defined
by the function $v_{0}$ (which can be thought of as a cdf). Then the measure 
$\lambda $ on $\boldsymbol{R}$ defined by the distribution function $%
v_{0}(F),$ namely $\lambda _{(v_{0}F)}\{(a,b]\}=v_{0}(F(b))-v_{0}(F(a)),$ is
the same as $\mu \varphi ^{-1}$ and the RHS of (\ref{32}) becomes $%
\int_{(-\infty ,a]}u_{0}(x)dv_{0}(F(x)).$ Furthermore the LHS becomes $%
\int_{(0,F(a)]}u_{0}(F^{-1}(\alpha ))dv_{0}(\alpha ).$ Hence 
\begin{equation*}
\int_{(-\infty ,a]}u_{0}(x)dv_{0}(F(x))=\int_{(0,F(a)]}u_{0}(F^{-1}(\alpha
))dv_{0}(\alpha )\qquad \text{for any given }a\in \boldsymbol{R.}
\end{equation*}%
A very similar proof yields:%
\begin{equation*}
\int_{(0,p)}v_{0}(\alpha )du_{0}(F^{-1}(\alpha ))=\int_{(-\infty
,F^{-1}(p))}v_{0}(F(x))du_{0}(x)\qquad \text{for any given }p\in (0,1]\text{ 
}.
\end{equation*}

Clearly we can replace $u_{0}(\cdot )$ by $u(\cdot )$ and $v_{0}(\cdot )$ by 
$v(\cdot ),$ and also let $a\rightarrow \infty $ and $p\rightarrow 1$, thus
we have%
\begin{equation}
\int_{-\infty }^{\infty }u_{0}(x)dv(F(x))=\int_{0}^{1}u_{0}(F^{-1}(\alpha
))dv(\alpha )  \label{CV1*}
\end{equation}

\begin{equation}
\int_{-\infty }^{\infty }u(x)dv_{0}(F(x))=\int_{0}^{1}u(F^{-1}(\alpha
))dv_{0}(\alpha )  \label{CV2*}
\end{equation}

\begin{equation}
\int_{-\infty }^{\infty }v_{0}(F(x))du(x)=\int_{0}^{1}v_{0}(\alpha
)du(F^{-1}(\alpha ))  \label{CV3*}
\end{equation}

\begin{equation}
\int_{-\infty }^{\infty }v(F(x))du_{0}(x)=\int_{0}^{1}v(\alpha
)du_{0}(F^{-1}(\alpha ))  \label{CV4*}
\end{equation}

\begin{lem}
For any cdf $F$ satisfying conditions $(a)$ of Theorem 1 we have 
\begin{equation*}
\int_{-\infty }^{x_{0}}v_{0}(F(x))du_{0}(x)+\int_{0}^{\alpha
_{0}}u_{0}(F^{-1}(\alpha ))dv_{0}(\alpha )=v_{0}(\alpha _{0})u_{0}(x_{0})
\end{equation*}%
whenever either $\alpha _{0}=F(x_{0})$ or $x_{0}=F^{-1}(\alpha _{0})$ (or
both). \newline
Furthermore, for all $F(x_{1}^{-})\leq \alpha _{1}\leq F(x_{1})$ and/or $%
F^{-1}(\alpha _{1})\leq x_{1}\leq F^{-1}(\alpha _{1}^{+})$ 
\begin{multline}
\int_{-\infty }^{x_{1}}v_{0}(F(x))du_{0}(x)+\int_{0}^{\alpha
_{1}}u_{0}(F^{-1}(\alpha ))dv_{0}(\alpha )=  \notag \\
=v_{0}(F(x_{1}))u_{0}(x_{1})+v_{0}(\alpha _{1})u_{0}(F^{-1}(\alpha
_{1}))-v_{0}(F(x_{1}))u_{0}(F^{-1}(\alpha _{1})).  \label{35}
\end{multline}
\end{lem}

\noindent Proof. All integrals are bounded because of $(a)$. We prove the
first statement. Fix $x_{0}$ and consider the case $\alpha _{0}=F(x_{0})$.
By change of variables in the second term and integration by parts 
\begin{eqnarray*}
\int_{-\infty }^{x_{0}}v_{0}(F(x))du_{0}+\int_{-\infty
}^{x_{0}}u_{0}(x)d(v_{0}F) &=&v_{0}(F(x_{0}))u_{0}(x_{0}) \\
&=&v_{0}(\alpha _{0})u_{0}(x_{0}).
\end{eqnarray*}%
Now fix $\alpha _{0}$ and assume $x_{0}=F^{-1}(\alpha _{0})$. Then again
applying change of variables and integration by parts but to the inverse $%
F^{-1}$ we have 
\begin{eqnarray*}
\int_{0}^{\alpha _{0}}v_{0}(a)d(u_{0}F^{-1})+\int_{0}^{\alpha
_{0}}u_{0}(F^{-1}(\alpha ))dv_{0} &=&v_{0}(\alpha _{0})u_{0}(F^{-1}(\alpha
_{0})) \\
&=&v_{0}(\alpha _{0})u_{0}(x_{0}).
\end{eqnarray*}%
To prove the second statement of the Lemma assume $F(x_{1}^{-})\leq \alpha
_{1}\leq F(x_{1}).$ Then $F^{-1}(\alpha _{1})=F^{-1}(F(x_{1}))$ and the
function $u_{0}(F^{-1}(\alpha ))$ is constant for all $\alpha _{1}\leq
\alpha \leq F(x_{1}).$ By the first part of the Lemma

\begin{center}
$%
\begin{array}{l}
\int_{-\infty }^{x_{1}}v_{0}(F(x))du_{0}(x)+\int_{0}^{\alpha
_{1}}u_{0}(F^{-1}(\alpha ))dv_{0}(\alpha )= \\ 
=v_{0}(F(x_{1}))u_{0}(x_{1})-\int_{\alpha
_{1}}^{F(x_{1})}u_{0}(F^{-1}(\alpha ))dv_{0} \\ 
=v_{0}(F(x_{1}))u_{0}(x_{1})-u_{0}(F^{-1}(\alpha
_{1}))[v_{0}(F(x_{1}))-v_{0}(\alpha _{1})] \\ 
=v_{0}(F(x_{1}))u_{0}(x_{1})+v_{0}(\alpha _{1})u_{0}(F^{-1}(\alpha
_{1}))-v_{0}(F(x_{1}))u_{0}(F^{-1}(\alpha _{1})).%
\end{array}%
$
\end{center}

Similarly when $F^{-1}(\alpha _{1})\leq x_{1}\leq F^{-1}(\alpha _{1}^{+})$,
which ends the proof of Lemma 4.\vspace{5mm}

\textbf{Remark \ }The above Lemma continues to hold replacing $v_{0}(\alpha
) $ by any $v(\alpha )\in V^{0}$ and/or $u_{0}(x)$ by any $u(x)\in U^{0}.$%
\vspace{5mm}

\noindent \textbf{Acknowledgments}. Much of this work was carried out while
the second author was Senior Visiting Fellow at the Istituto di Studi
Avanzati of the University of Bologna in 2007 and during subsequent visits
to the Department of Statistical Sciences of Bologna University. The first
author wishes to thank the hospitality of the Department of Statistics at
the London School of Economics during a number of visits. Both authors would
like to thank Irene Crimaldi for providing the elegant proof of Theorem 4,
Claudio Zoli for useful conversations at an early stage of the work and two
referees for helpful comments and additional references.\vspace{3mm}

\begin{center}
{\Large References}
\end{center}

Atkinson, A.B. (1970). On the measurement of inequality. \textit{Journal of
Economic Theory}, \textbf{2}, 244-263.

Arrow, K.J. (1974). \textit{Essays in the Theory of Risk-Bearing \ }%
North-Holland: New York

Cowell, F.A. (1977). \textit{Measuring Inequality}. Oxford: Philip Allan.

Chateauneuf, A. and Moyes, P. (2005). Does the Lorenz curve really measure
inequality? Another look at inequality measurement. Econometric Society
World Congress, University College London, London, UK, August 2005.

Donaldson, D. and Weymark, J. (1980). A single-parameter generalization of
Gini indices of inequality, \textit{Journal of Economic Theory}, \textbf{2},
67-86.

Lorenz, M.O. (1905). Methods of measuring concentration of wealth. \textit{%
Journal of\ the American Statistical Association,} \textbf{9}, 209--219.

Kahneman, D. and Tversky, A. (1979). Prospect theory: an analysis of
decision under risk. \textit{Econometrica}, \textbf{47}, 263-291.

Kahneman, D. and Tversky, A. (1992). Advances in prospect theory; cumulative
representation of uncertainty. \textit{Journal of Risk Uncertainty}, \textbf{%
5}, 297-323

Maccheroni, F., Muliere, P and Zoli, C. (2005). Inverse stochastic orders
and generalized Gini functionals. \textit{Metron}, \textbf{63}, 529-559.

Marshall, A.W. and Olkin, I. (1979). \textit{Inequalities: Theory of
Majorization and its Applications}. New York: Academic Press.

Muliere, P. and Scarsini, M. (1989). A note on stochastic dominance and
inequality measures. \textit{Journal of Economic Theory}, \textbf{49},
314-323.

M\"{u}ller, A. and Stoyan, D. (2002). \textit{Comparison Methods for
Stochastic Models and Risks}. New York: Wiley.

Ogryczak, W. and Ruszczynski, A. (2002). Dual stochastic dominance and
related mean--risk models. \textit{SIAM Journal of Optimization}, \textbf{13}%
, 60-78.

Pallaschke, D. and Rolewicz, S. (1997).\textit{Foundations of Mathematical
Optimisation: Convex Analysis without Linearity}. Kluwer: Dordrecht.

Pratt, J.W. (1964).\ Risk aversion in the small and in the large. \textit{%
Econometrica,}\textbf{\ 32}, 122-136

Quiggin, J. (1982). A theory of anticipated utility. \textit{Journal of
Economic Behaviour and Organization}, \textbf{3}, 323--343.

Sen, A.K. (1973). \textit{On Income Inequality}. Oxford: University Press.

Shaked, M. and Shanthikumar, J.G. (2007). \textit{Stochastic Orders}. New
York: Academic Press.

Sordo, M.A. and Ramos, H.M. (2007). Characterization of stochastic orders by
L-functionals. \textit{Statistical Papers}, 48, 249-263

Yaari, M.E. (1987). The dual theory of choice under risk. \textit{%
Econometrica}, \textbf{55}, 95--115.

Yitzhaki, S. (1983). On an extension of the Gini inequality index. \textit{%
International Economic Review}, \textbf{24}, 617-628.

\end{document}